\documentclass[12pt,leqno]{amsart}
\usepackage{amssymb,verbatim,enumerate,ifthen,hyperref}
\usepackage{mathtools} 

\usepackage{cite}
\usepackage{tikz-network}
\usepackage{graphicx}
\usepackage[mathscr]{eucal}
\usepackage[utf8]{inputenc}
\usepackage[T1]{fontenc}
\def\N{\mathbb{N}}
\def\R{\mathbb{R}}

\def\C{\mathscr{C}}

\long\def\comment#1{}

\newtheorem{theorem}{Theorem}[section]
\newtheorem*{theorem*}{Theorem}
\def\Thm#1#2{\ifthenelse{\equal{#1}{*}}{\begin{theorem*}#2\end{theorem*}}
             {\begin{theorem}\label{T#1}#2\end{theorem}}}
\newtheorem{Atheorem}{Theorem}

\newtheorem{proposition}[theorem]{Proposition}
\newtheorem*{proposition*}{Proposition}
\def\Prp#1#2{\ifthenelse{\equal{#1}{*}}{\begin{proposition*}#2\end{proposition*}}
{\begin{proposition}\label{P#1}#2\end{proposition}}}
\def\prp#1{Proposition~\ref{P#1}}

\newtheorem{corollary}[theorem]{Corollary}
\newtheorem*{corollary*}{Corollary}
\def\Cor#1#2{\ifthenelse{\equal{#1}{*}}{\begin{corollary*}#2\end{corollary*}}
             {\begin{corollary}\label{C#1}#2\end{corollary}}}
\def\cor#1{Corollary~\ref{C#1}}

\newtheorem{lemma}[theorem]{Lemma}
\newtheorem*{lemma*}{Lemma}
\def\Lem#1#2{\ifthenelse{\equal{#1}{*}}{\begin{lemma*}#2\end{lemma*}}
             {\begin{lemma}\label{L#1}#2\end{lemma}}}

\theoremstyle{definition}
\newtheorem{remark}[theorem]{Remark}
\newtheorem*{remark*}{Remark}
\def\Rem#1#2{\ifthenelse{\equal{#1}{*}}{\begin{remark}\rm #2\end{remark}}
             {\begin{remark}\label{R#1}\rm #2\end{remark}}}

\newtheorem{example}[theorem]{Example}
\newtheorem*{example*}{Example}
\def\Exa#1#2{\ifthenelse{\equal{#1}{*}}{\begin{example*}\rm #2\end{example*}}
             {\begin{example}\label{Ex#1}\rm #2\end{example}}}

\def\eq#1{{\rm(\ref{E#1})}}
\def\Eq#1#2{\ifthenelse{\equal{#1}{*}}
  {\begin{equation*}\begin{aligned}#2\end{aligned}\end{equation*}}
  {\begin{equation}\begin{aligned}\label{E#1}#2\end{aligned}\end{equation}}}

\begin{document}
\begin{flushright}
\end{flushright}
\vspace{5mm}

\date{\today}

\title[On Weighted Star--Convex Graphs]
{On Weighted Star--Convex Graphs}

\author[A. R. Goswami]{Angshuman R. Goswami}
\address[A. R. Goswami]{Department of Mathematics, University of Pannonia,
H-8200 Veszprem, Hungary}
\email{goswami.angshuman.robin@mik.uni-pannon.hu}


\thanks{The research of the author was supported by the 
``Ipar a Veszprémi Mérnökképzésért'' Foundation
}

\begin{abstract}
The primary objective of this paper is to investigate the notions of geometric and sequential convexity within a graph-theoretic framework, with the aim of examining various structural properties and exploring the connection between these two branches of mathematics.\\

A simple connected vertex-weighted graph $G(V,E)$ with a non-empty set of leaf vertices is said to be star-convex if there exists at least one node $u\in V(G)$ such that, for every chosen leaf vertex $v$, there is a monotone path (either increasing or decreasing) connecting $v$ to $u$. One of the main results states that a graph $G$ is star-convex if and only if there exists a tree $T\subseteq G$ that contains all leaf vertices and is itself star-convex.\\

On the other hand, a sequence $\big(u_n\big)_{n=0}^{\infty}$ is said to be convex if it satisfies the following inequality
\Eq{*}{
2u_{i}\leq u_{i-1}+u_{i+1}\qquad \mbox{for all}\quad i\in\N.
}
We demonstrate that, under minimal assumptions, a class of convex sequences can be embedded into a spider graph so as to make it star-convex.\\

The necessary definitions, background, motivation, and other relevant details are discussed in the following section.
\end{abstract}

\maketitle
\section*{Introduction}
Throughout this paper, $\mathbb{R}_{+}$ denotes the set of non-negative real numbers. 
Let $G(V,E)$ be a simple, connected, unweighted graph with vertex set $V$ and edge set $E$. If a weight function $w:V\to \mathbb{R}_{+}$ is assigned, then the resulting weighted graph is denoted by $G(V,E,w)$. It is important to note that the graphs $G(V,E,w_{1})$ and $G(V,E,w_{2})$ have the same vertex and edge sets but differ in their corresponding weight functions. Furthermore, $V_{\ell}(\subseteq V)$ denotes the set of all vertices of $G$ with degree exactly $1$. Similarly, a tree is denoted by $T(V,E)$, and the associated graph-theoretic terminology can be expressed using analogous notation.\\

In convex geometry, a set $X$ is said to be \textit{star-convex} if there exists a point $p\in X$ such that, for every $x\in X$, the unique line segment connecting $p$ and $x$ lies entirely within $X$. The subset consisting of all points having the same property as $p$ is called the \textit{center of the set} $X$. It is clear that every convex set is star-convex; however, the converse does not hold in general. A basic example of it is a connected compact star-shaped set. The study of star-convex bodies is important in optimization and related algorithmic questions. Further background in this direction can be found in \cite{KRASNOSEL,Krein,Stanek,Breen,Hansen} and the references therein.\\

In graph theory, a tree $T(V,E)$ having exactly one internal node $v$ of degree at least $3$ may be regarded as the graph-theoretic analogue of geometric or set-theoretic star-convexity. In the literature, this type of tree is usually referred to as a spider graph. In other words, a \textit{spider} is a tree with exactly one vertex of degree greater than $2$, while all other internal vertices have degree exactly $2$. In a spider graph, the vertex of largest degree is called the \textit{hub}, and each maximal path from the hub to a leaf is called a \textit{leg}. Thus, a star is the special case in which every leg has length exactly $1$, so that a spider is a direct generalization of a star and a proper subclass of trees. If all legs of a spider have equal length, then the spider is called a \textit{balanced spider or regular spider}. A number of structural and extremal results on spider graphs, including work related to the Erd\H{o}s--S\'os conjecture, may be found in \cite{Fan,Gargano,Jiang,Yang,Bonato}.\\

On the other hand, for any non-empty and non-singleton convex set $I\subseteq\R$, a function $f:I\to\R$ is said to be \textit{star-convex} if there exists a point $p\in I$ such that, for every $x\in I$, the line segment connecting $(p,f(p))$ and $(x,f(x))$ lies entirely either in the epigraph of $f$ or in the hypograph of $f$. The notion of a star-convex function not only generalizes ordinary convexity but also highlights its relationship with star-convex sets. For further details, we refer to \cite{Pales,Hinder,Hansen}.\\
 
These interconnected concepts motivate us to analyse graphs that possess a spider-like structure by assigning weights to their vertices in a manner resembling a discrete convex or star-convex function. In other words, we aim to develop a graph-theoretic framework that mirrors the well-established features of star-convexity in function theory. We first introduce the notion of star-convex graphs and then discuss the concept of discrete convexity in order to build a bridge between these two settings.\\

A weighted graph $G(V,E,w)$ with non-empty $V_{\ell}$ is said to be \textit{star-convex} if there exists at least one vertex $u\in V$ such that, for every $v\in V_{\ell}$, there exists a weighted-monotone path $P$ connecting $u$ and $v_{\ell}$. In other words, if $P:=\{u,v_1,\cdots,v_{\ell}\}$, then one of the following holds:
\Eq{*}{
w(u)\leq w(u_1)\leq \cdots\leq w(v_{\ell})\qquad\mbox{or}\qquad 
w(u)\geq w(u_1)\geq \cdots\geq w(v_{\ell}).
}
The collection of all such vertices $u$ is called the \textit{core} of the graph $G$, and we denote it by $C(G)$. Thus, the core consists of those vertices from which a weighted-monotone path can be traced to every vertex of degree $1$. Clearly, all the weighted graphs with exactly one leaf vertex turns to be star-convex. Hence in our research, we are interested to investigate only those star-convex graphs which contain atleast two leaf nodes.\\

Besides, the motivation for studying weighted star-convexity is even more natural. For example, regardless of how weights are assigned to the vertices of a star graph, the resulting weighted graph is always star-convex. We therefore investigate several structural properties of weighted graphs possessing star-convexity. In particular, we show that both the union and the intersection of two graphs $G_1$ and $G_2$ can be star-convex provided that their cores have a non-empty intersection. We study star-convexity on trees and show that the extremal weights are attained either at leaf vertices or at vertices from the core. As one of our main results, we prove that a weighted graph $G(V,E,w)$ is star-convex if and only if there exists a star-convex tree $T(V_{\ell},E,w)$ extractable from $G$.\\

One of the key objectives of this paper is to relate two distinct objects from different branches of discrete mathematics, convex sequences and spider graphs. A real-valued sequence $\big(u_n\big)_{n=0}^\infty$ is called \textit{convex} if the following discrete functional inequality is satisfied
\Eq{*}{
u_{i+1}-u_{i}\leq u_{j+1}-u_{j}\qquad \mbox{for all $i,j\in\N$}\quad\mbox{with}\quad i\leq j.
}
The study of sequential convexity is important because every convex sequence can be viewed as a discrete analogue of a convex function, which makes it relevant in several discrete optimization problems. Moreover, many well-known sequences, such as arithmetic, geometric, partition-related, and Fibonacci-type sequences, may be studied within the broader framework of convex sequences. Various generalized and approximate forms of convex sequences have also been investigated, leading to several notable characterizations, structural properties, and algorithmic applications; see, for example, \cite{Mitrinovicc,Essen,pecaric,Debnath,GauSte} and the references therein.\\

We show that it is possible to embed a class of convex sequences into a spider graph of suitable cardinality to obtain a star-convex graph. This result may also be interpreted as a transportation problem for weight allocation.\\

We start our investigation with some basic structural properties.

\section{Main Results}
The following proposition is straightforward, and therefore, we record only its statement.
\Prp{0}{Let $G_1$ and $G_2$ be two star-convex graphs. If $G_1\subset G_2$ then $c(G_1)\subseteq c(G_2)$.
}
In most cases, without any additional assumptions, the intersection of two algebraic structures or two geometric figures of the same type tends to preserve their essential behavioral properties. In contrast, their union usually fails to retain those core properties. This is why one being contained within the other is often a necessary condition for the union to behave well. This pattern holds for structures such as groups, rings, convex bodies, and so on. However, the results established in the following proposition go against this usual norm.\\
\Prp{1}{
Let $G_1$ and $G_2$ be two star-convex graphs such that both vertex sets $V_{_{\ell}}(G_1\cap G_2)$ and $V_{_{\ell}}(G_1\cup G_2)$ are non-empty. Then the following two assertions can be made
\begin{enumerate}[(i)]
\item If $C(G_1)\cap C(G_2)$ is non-empty then $G_1\cup G_2$ is star-convex.
\item If $C(G_1)\cap C(G_2)$ is non-empty then $G_1\cap G_2$ may not be star-convex.
\end{enumerate}
}
\begin{proof}
To prove the first assertion, we assume $u\in C(G_1)\cap C(G_2)$ arbitrarily. This implies $u\in C(G_1)$ and $u\in C(G_2)$. Now, let $v\in V_{_{\ell}}(G_1\cup G_2)$, which indicates either $v\in V_{_{\ell}}(G_1)$ or $v\in V_{_{\ell}}(G_2)$. Without loss of generality, we consider $v\in V_{_{\ell}}(G_1)$. Thus, there exists a weighted-monotone path $P$ connecting the vertices $u$ and $v$. Since $u\in C(G_1)\cap C(G_2)$ and $v\in V_{_{\ell}}(G_1\cup G_2)$ are arbitrary, we conclude that $G_1\cup G_2$ possesses star-convexity.\\

The second assertion of the above proposition is illustrated with the following example.
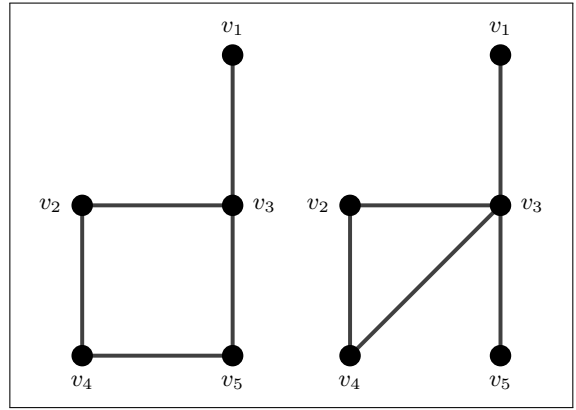
\begin{center}
\fbox{
\begin{tikzpicture}
\Vertex[x=0,y=0, size=0.25, color=black, label=$v_1$,  position=above]{A}
\Vertex[x=-2,y=-2, size=0.25, color=black, label=$v_2$,  position=left]{B}
\Vertex[x=0, y=-2, size=0.25, color=black, label=$v_3$,  position=right]{C}
\Vertex[x=-2, y=-4,  size=0.25, color=black, label=$v_4$,  position=below]{D}
\Vertex[x=0, y=-4, size=0.25, color=black, label=$v_5$,  position=below]{E}
\Edge(A)(C)
\Edge(B)(C)
\Edge(B)(D)
\Edge(C)(E)
\Edge(D)(E)
\end{tikzpicture}

\begin{tikzpicture}
\Vertex[x=24,y=0, size=0.25, color=black, label=$v_1$,  position=above]{A}
\Vertex[x=22,y=-2, size=0.25, color=black, label=$v_2$,  position=left]{B}
\Vertex[x=24, y=-2, size=0.25, color=black, label=$v_3$,  position=right]{C}
\Vertex[x=22, y=-4,  size=0.25, color=black, label=$v_4$,  position=below]{D}
\Vertex[x=24, y=-4, size=0.25, color=black, label=$v_5$,  position=below]{E}
\Edge(A)(C)
\Edge(B)(C)
\Edge(B)(D)
\Edge(C)(E)
\Edge(C)(D)
\end{tikzpicture}
}
\end{center}
\begin{center}
$w(v_1)=1$, $w(v_2)=1$, $w(v_3)=2$, $w(v_4)=2$, $w(v_5)=2$ 
\\
\end{center}
The two graphs above share the same vertex sets and the corresponding vertex weights. The two graphs preserve star convexities with the core $v_3$. However, one can easily check that the intersection of these two graphs fails to meet the criteria of star-convexity.
\end{proof}

The following two results show that if a leaf vertex of a star-convex tree also belongs to its core, then the extremal weights are attained at leaf vertices.
\Prp{2}{Let $T(V,E,w)$ be a star-convex tree such that $u\in C(T)\cap V_{_{\ell}}(T)$. Then for all $v\in V_{_{\ell}}$, the path connecting $u$ and $v$ has the same weighted monotonicity. 
}
\begin{proof}
To prove the result, let $P=\{u,v_1,\cdots,v_n,v\}$ be the path connecting the distinct vertices $u$ and $v$. Since $u\in C(G)$, the path $P$ is weighted-monotone. Without loss of generality, we assume that this monotonicity is increasing. That is, $P$ satisfies the following inequality:
\Eq{101}{
w(u)\leq w(v_1)\leq \cdots \leq w(v_n)\leq w(v).
}
Now consider another distinct vertex $v'$, and let the corresponding path connecting $u$ and $v'$ be denoted by $P'=\{u,v'_1,\cdots,v'_m,v'\}$. Since $u\in V_{_{\ell}}(T)$, the paths $P$ and $P'$ imply that $v_1=v'_1$. Hence, the condition $u\in C(T)$ together with \eq{101} yields the following
\Eq{*}{
w(u)\leq w(v_1)=w(v'_1)\leq \cdots \leq w(v'_m)\leq w(v').
}
Thus, the same assertion holds for all $v\in V_{_{\ell}}$. This completes the proof.
\end{proof}
\Cor{3}{Let $T(V,E,w)$ be a star-convex tree such that $C(T)\cap V_{_{\ell}}(T)$ is non-empty. Then the following holds
\Eq{105}{
\max_{v\in V}\Big(w(v)\Big)=\max_{v\in V_{_{\ell}}\cup\, C}\Big(w(v)\Big)\qquad\mbox{and}\qquad
\min_{v\in V}\Big(w(v)\Big)=\min_{v\in V_{_{\ell}}\cup \,C}\Big(w(v)\Big).
}
}
\begin{proof}
Let $u\in C(T)\cap V_{_{\ell}}(T)$ be arbitrary. Then using \prp{2} and without loss of generality, first we assume that for any distinct $v\in V_{_{\ell}}$, the path $P$ connecting the vertices $u$ and $v$ satisfies \eq{101}. Also, since a tree is a simple connected graph without cycles, the union of all such paths $P$ traces all the internal vertices of $T$. Hence, from the monotone increasing property of all these paths, we conclude that
\Eq{102}{
\max_{v\in V}\Big(w(v)\Big)=\max_{v\in V_{_{\ell}}\setminus{u}}\Big(w(v)\Big)\qquad\mbox{and}\qquad
\min_{v\in V}\Big(w(v)\Big)=\min_{v\in V_{_{\ell}}}\Big(w(v)\Big)=w(u).
}
In the case of the reverse inequality of \eq{101}, we obtain 
\Eq{103}{
\max_{v\in V}\Big(w(v)\Big)=\max_{v\in V_{_{\ell}}}\Big(w(v)\Big)=w(u)\qquad\mbox{and}\qquad
\min_{v\in V}\Big(w(v)\Big)=\min_{v\in V_{_{\ell}}\setminus{u}}\Big(w(v)\Big).
}
The relations \eq{102} and \eq{103} cover all possibilities. This validates \eq{105} and completes the proof.
\end{proof}
In fact, using the above corollary, we can propose the following generalized result.
\Prp{4410}{Let $T(V,E,w)$ be a star-convex tree. Then the following equalities holds
\Eq{*}{
\max_{v\in V}\Big(w(v)\Big)=\max_{v\in V_{_{\ell}}\cup\, C}\Big(w(v)\Big)\qquad\mbox{and}\qquad
\min_{v\in V}\Big(w(v)\Big)=\min_{v\in V_{_{\ell}}\cup \,C}\Big(w(v)\Big).
}
}
\begin{proof}
If the vertex set $C(T)\cap\big(V(T)\setminus V_{_{\ell}}(T)\big)$ is non-empty, one can also easily verifies the validity of \eq{105}. This together with \cor{3} establishes this assertion.
\end{proof}
The next proposition shows that the structural analysis of a star-convex graph can be reduced to the study of a tree extracted from it, which contains all of its leaf vertices and traces its core.
\Thm{201}{ The graph $G(V,E,w)$ is star-convex if and only if there exists a star-convex tree $T\subseteq G$ such that $V_{_{\ell}}(G)=V_{_{\ell}}(T)$ holds.
}
\begin{proof}
First we assume, $T\subseteq G$ is a star-convex tree satisfying the condition $V_{_{\ell}}(G)=V_{_{\ell}}(T)$. The star-convexity of $T$ implies $C(T)$ is non-empty, and hence we assume $u\in C(T)$. Thus, we can conclude that for every  $v\in V_{_{\ell}}(G)$ there exists a monotone path $P$ that connects $u$ with $v$. This shows that the graph $G(V,E,w)$ possesses star-convexity.\\

For the converse, assume that $G$ is star-convex. Then $C(G)$ is non-empty, and let $u\in C(G)$. If $G$ is a tree, the result is immediate. Hence, we assume that $G$ is not a tree. For each $v_{i}\in V_{_{\ell}}(G)$, choose a monotone weighted path $P_i$ connecting $u$ and $v_i$. If the union of all such paths $\cup{P_i}$ is a tree, then the result follows. Otherwise, $\cup{P_i}\subseteq G$ must contain a cycle. Without loss of generality, assume that there exist $v_i,v_j\in V_{_{\ell}}(G)$ with corresponding monotone paths $P_i$ and $P_j$ such that $P_i\cup P_j$ contains a cycle. In order for a cycle to appear, the paths $P_i$ and $P_j$ must have at least two vertices in common, one of which is $u$. Let the other common vertex be $u'$. For simplicity, we sketch the two weighted-monotone paths $P_i$ and $P_j$ as follows:
$$
P_i=u--u'--v_i\qquad \mbox{and}\qquad P_j=u--u'--v_j.
$$
Since the paths meet at $u'$, they must have the same monotonicity along the common segment. Therefore, instead of following two separate routes, the paths $P_i$ and $P_j$ may be modified to follow the same route from $u$ to $u'$. This removes the cycle without affecting our objective. Repeating this procedure eliminates all cycles and leaves a tree. Clearly, in the resulting tree, $u\in C(T)$, and the condition $V_{_{\ell}}(G)=V_{_{\ell}}(T)$ remains satisfied. This completes the proof. 
\end{proof}
In the next result, we establish the relationship between a regular spider graph and the class of ordinary convex sequences. The same idea can be adapted to more general spider graphs by adjusting the cardinalities of the convex sequences appropriately. One may also attempt to extend this connection to more general graph classes such as lobster graphs or scorpion graphs. It is also important to mention that for any sequence $\big(u_n\big)_{n=0}^{\infty}$, the symbol $\{u\}$ represents its range set.
\Thm{1000}{
Let $\C:=\bigg\{\Big(u^{^{j}}_{_{i}}\Big)_{i=1}^{2\ell+1}\,\,\Big|\,\,j=1,\cdots,n\bigg\}$ be the class of $n$-convex sequences satisfying the following two conditions
\Eq{1005}{
(i)\,\,\,u^{^{j}}_{_{\ell+1}}=u \quad\mbox{for all}\quad j\in\{1,\cdots,n\}\qquad \mbox{and}\qquad (ii)\,\,\, \min_{j\in[1,n]\cap\N}\Big\{u^{^{j}}\Big\}={u}.
}
We assume $S(V,E)$ is a regular spider graph with $2n$ legs and the length of each leg is $\ell$. Then the points from the sequences in $\C$ can be embedded into $V(S)$ as weights to make it a star-convex graph.
}
\begin{proof}
Let $v\in V$ be the hub vertex of the spider graph $S(V,E)$.
To prove the theorem, we consider $v_1$ and $v_2$ are two distinct leaf vertices and the path $P_1$ that connects these two vertices, given by $P_1:=\{v_1=v_{_{1}}^{^{1}},v_{_{2}}^{^{1}}\cdots,v_{_{2\ell+1}=v_2}^{^{1}}\}$.\\

Since a spider graph is a special type of tree, this implies the path $P_1$ is unique. Besides, as per the property of the spider, it has to pass through the hub vertex $v$. Now, without loss of generality, we choose the sequence $\Big(u^{^{1}}_{_{i}}\Big)_{i=1}^{2\ell+1}$ and define the mapping $w_1:P_1 \to \Big\{u^{^{1}}\Big\}$ as follows
$$w_1\Big(v_{_{i}}^{^{1}}\Big)=u_{_{i}}^{^{1}}\qquad\mbox{for all}\quad i\in[1,2\ell+1]\cap\N\,.$$
Next, we choose another pair of distinct leaf vertices (say $v_3$ and $v_4$) and the path $P_2$ connecting those two nodes. The hub vertex $v$ is the only vertex that the two paths $P_1$ and $P_2$ have in common. Besides, it is the 
$(\ell+1)^{th}$ vertex in both paths. Therefore, using the same framework, we can select the sequence 
$\Big(u^{^{2}}_{_{i}}\Big)_{i=1}^{2\ell+1}$ and able to define the weight function $w_2:P_2 \to \Big\{u^{^{2}}\Big\}$ as  \eq{1005} ensures the equality $w_{_{1}}\Big(v_{_{\ell+1}}^{^{1}}\Big)=w_{_{2}}\Big(v_{_{\ell+1}}^{^{2}}\Big)=w(v)=u$. That is no overlapping of assigned weights. Continuing this methodology, we allocate weights to each and every vertex of S and eventually arrive at the weighted spider graph 
$S(V,E,w)$.\\

Now choose an arbitrary leaf vertex $v_k$ and consider the path $P$ connecting the root $u$ and $v_k$. To establish the star-convexity of $S(V,E,w)$, it is sufficient to show that the weighted path $P$ is monotone. By construction, the vertex $v_k$ is assigned the extremal sequential value of one of the sequences in $\C$. Let this sequence be $\Big(u^{^{j}}_{_{i}}\Big)_{i=1}^{2\ell+1}$ for some $j\in\{1,\cdots,n\}$.\\

Using the sequential convexity property and the two conditions stated in \eq{1005}, we obtain
\Eq{*}{
u^{^{j}}_{_{2}}-u^{^{j}}_{_{1}}\leq\cdots\leq u^{^{j}}_{_{n}}-u^{^{j}}_{_{n-1}}\leq 0\leq u^{^{j}}_{_{n+1}}-u^{^{j}}_{_{n}}\leq \cdots\leq u^{^{j}}_{_{2\ell+1}}-u^{^{j}}_{_{2\ell}}
\,.} 
From the inequality above, we get the following
\Eq{*}{
u^{^{j}}_{_{n}}\leq\cdots\leq u^{^{j}}_{_{2}}\leq u^{^{j}}_{_{1}}
\qquad \mbox{and}\qquad
u^{^{j}}_{_{n}}\leq u^{^{j}}_{_{n+1}}\leq \cdots \leq u^{^{j}}_{_{2\ell+1}}.
}
In other words, the weighted path $P$ is either monotonically increasing or monotonically decreasing. This yields the star-convexity of 
$S(V,E,w)$ and completes the proof.
\end{proof}
This research offers numerous avenues for further development, with potential applications in fields such as chemical graph theory, network science, decision theory, and algorithmic studies. For example, many molecules exhibit structures resembling spider graphs, where star convex like vertex-weighted analysis could prove useful in understanding their chemical reactions with other compounds. Similarly, the design of algorithms for detecting the core of a star-convex graph presents an interesting as well as challenging direction.
\section*{Statements and Declarations}

\noindent\textbf{Funding.} 
The author received financial support from the \textit{``Ipar a Veszprémi Mérnökképzésért'' Foundation} for the research presented in this work.

\vspace{6pt}
\noindent\textbf{Competing Interests.} 
The author declare that there are no financial or non-financial competing interests relevant to the contents of this article.

\vspace{6pt}
\noindent\textbf{Ethics Approval.} 
Not applicable. This study does not involve human participants or animals.

\vspace{6pt}
\noindent\textbf{Consent to Participate.} 
Not applicable.

\vspace{6pt}
\noindent\textbf{Consent for Publication.} 
Not applicable.

\vspace{6pt}
\noindent\textbf{Data, Materials and/or Code Availability.} 
No datasets or code were generated or analysed during the current study.

\end{document}